\documentclass[11pt]{article}
\usepackage{amsmath,amssymb,amsthm,amsfonts}
\numberwithin{equation}{section}
\usepackage{mathtools,mathrsfs,bm,graphicx,float}
\usepackage{multirow,booktabs,longtable,arydshln,enumerate,cite,microtype,enumitem}
\usepackage{hyperref}
\usepackage[a4paper,margin=2.5cm]{geometry}
\hypersetup{colorlinks=true,linkcolor=blue,citecolor=blue,urlcolor=blue,
 pdftitle={Weak Approximate Spectral Synthesis in Some Analytic Hilbert Spaces},
 pdfauthor={},pdfsubject={Kernel comparison, bounded coordinate multipliers and finite spectral approximation}}
\allowdisplaybreaks[2]
\newcommand{\doi}[1]{\href{https://doi.org/#1}{\nolinkurl{doi:#1}}}
\newtheorem{theorem}{Theorem}[section]
\newtheorem{proposition}[theorem]{Proposition}
\newtheorem{lemma}[theorem]{Lemma}
\newtheorem{corollary}[theorem]{Corollary}
\newtheorem{remark}[theorem]{Remark}
\newcommand{\D}{\mathbb D}
\newcommand{\norm}[1]{\left\lVert#1\right\rVert}
\newcommand{\ip}[2]{\langle#1,#2\rangle}
\title{\Large\bfseries
Weak Approximate Spectral Synthesis in Radial Weighted Bergman Spaces and Hyponormal  Weighted Shifts}

\author{Zhaopeng Lin\thanks{Corresponding author.}, Yufeng Lu, Chao Zu }
\date{}

\begin{document}
\maketitle
\begin{center}
\begin{minipage}{14cm}\small
\noindent\textbf{Abstract}\quad 
We establish weak approximate spectral synthesis
for reproducing kernel Hilbert spaces on the unit disc
with radial kernels for which the coordinate multiplier
is a hyponormal contraction. We further extend this result to reproducing
kernel Hilbert spaces whose kernels are
comparable to such radial kernels in the
positive semidefinite order. 
As an application, we answer Shimorin's question
[Duke Math. J. 101 (2000), 1--39, Section~5]concerning weak approximate spectral synthesis
in general radial weighted Bergman spaces. 
In particular, we extend his result
[Proc. Amer. Math. Soc. 131 (2003), 1777--1787]
for standard weighted Bergman spaces
$A_\alpha^2$ from $-1<\alpha\le1$
to the full range $\alpha>-1$.

\medskip

\noindent\textbf{Mathematics Subject Classification (2020):}
Primary 47A15, 46E22;
Secondary 30H20, 47B32, 47B37.

\noindent\textbf{Keywords.}\quad
Weak approximate spectral synthesis,
hyponormal weighted shift,
reproducing kernel,
weighted Bergman space.
\end{minipage}
\end{center}

\section{Introduction}\label{sec:intro}

The study of invariant subspaces of the shift
operator is an important topic in operator
theory and analytic function theory. 
Let $\mathcal H(K)\subset\operatorname{Hol}(\D)$
be a reproducing kernel Hilbert space on the
unit disc $\D$. The shift operator $S_{\mathcal H(K)}$
is defined by
\[
(S_{\mathcal H(K)}f)(z)=zf(z),
\qquad f\in\mathcal H(K).
\]
On the Hardy space $H^2$, Beurling's theorem
completely characterizes the
$S_{H^2}$-invariant subspaces: every nonzero
closed $S_{H^2}$-invariant subspace is of
the form $\theta H^2$ for some inner
function $\theta$.
For Bergman and other analytic Hilbert spaces,
however, the structure of shift-invariant
subspaces is more complicated
\cite{ARS,Shimorin2000}.

Since a closed subspace $M$ is
$S_{\mathcal H(K)}$-invariant if and only if
$M^\perp$ is $S_{\mathcal H(K)}^*$-invariant,
we consider the approximation problem for
$S_{\mathcal H(K)}^*$.
A natural class of finite-dimensional
$S_{\mathcal H(K)}^*$-invariant subspaces
arises from finite zero sets
\cite{Hedenmalm1991,DKSS1993}.
For a finite set $A\subset\D$, consider
\[
I_A=\{f\in\mathcal H(K):
f(\lambda)=0,\ \lambda\in A\}.
\]
Then $I_A$ is a closed
$S_{\mathcal H(K)}$-invariant subspace.
By the reproducing property,
\[
I_A^\perp
=\operatorname{span}\{K_\lambda:\lambda\in A\},
\]
where $K_\lambda$ denotes the reproducing
kernel of $\mathcal H(K)$ at $\lambda\in\D$.
Moreover,
\[
S_{\mathcal H(K)}^*K_\lambda
=\overline{\lambda}K_\lambda,
\qquad \lambda\in\D.
\]
Thus $I_A^\perp$ is a finite-dimensional
$S_{\mathcal H(K)}^*$-invariant subspace
spanned by eigenvectors.

More generally, a nonzero vector $v$ is called
a root vector of an operator $T$ if
\[
(T-\zeta I)^m v=0
\]
for some $\zeta\in\mathbb C$ and integer
$m\ge1$. Eigenvectors correspond to $m=1$.  
Spectral synthesis asks whether every
closed subspace invariant under $T$ can be recovered
as the closed linear span of the root vectors
it contains. As discussed by Nikol'skii
\cite[\S1]{Nikolskii}, this property fails
even for the backward shift on $H^2$.
Indeed, if $\theta$ is a nonconstant singular
inner function, then $
H^2\ominus\theta H^2$ 
is a nonzero subspace invariant under the backward
shift and contains no nonzero root vectors.

Motivated by this failure, Nikol'skii
\cite[\S2]{Nikolskii} proposed in 1978
a weaker approximation problem:
 can every closed
subspace invariant under $S_{\mathcal H(K)}^*$
be obtained as a limit of invariant subspaces
of finite dimension, without requiring
the approximating subspaces to be increasing
or contained in the subspace?

More precisely, we consider the formulation of Nikol'skii's question
in a Hilbert space. Given a closed
subspace $N\subset\mathcal H(K)$ invariant
under $S_{\mathcal H(K)}^*$, does there exist a sequence
of subspaces $E_j$ of finite dimension invariant
under $S_{\mathcal H(K)}^*$ such that
\[
N=\liminf_{j\to\infty}E_j,
\]
where
\begin{equation}\label{eq:lowerlimit}
\liminf_{j\to\infty}E_j
=
\bigl\{x\in\mathcal H(K):
\operatorname{dist}_{\mathcal H(K)}(x,E_j)
\longrightarrow0\bigr\}?
\end{equation}
Following Shimorin \cite[Definitions~1--4]{Shimorin2000},
this property is called weak approximate spectral
synthesis.

For the Hardy space $H^2$, weak approximate
spectral synthesis follows from the work of
Douglas, Shapiro and Shields \cite{DSS};
see also \cite[p.~37]{Nikolskii1986}.
Shimorin \cite{Shimorin2000} proved that
weak approximate spectral synthesis also holds
for $A^2$.
Subsequently, Shimorin \cite{Shimorin2003}
extended this property to certain weighted
Bergman spaces. More precisely, the property
holds for radial weights $\omega$ satisfying
\[
\Delta\log\frac{\omega(z)}{1-|z|^2}\ge0.
\]

In particular, this condition is satisfied
by the standard weighted Bergman spaces
$A_\alpha^2$ for $-1<\alpha\le1$.
However, the case $\alpha>1$ remained open. 
More generally, as Shimorin observed in
\cite[\S5]{Shimorin2000}, his method cannot
be directly extended to arbitrary radial
weighted Bergman spaces, since the required
positivity of the associated weighted
biharmonic Green function may fail.

The Hardy operator is closely related to the
backward shift on $H^2$ \cite{BHS1965},
and its invariant subspaces have been studied
through finite-dimensional approximation
\cite{AglerMcCarthy2023}. 
Fu, Guo and Yan \cite{FGY} established
coanalytic models relating the logarithmic
Hardy integral operator $H_\beta$ on
$L^2[0,1]$, defined by
\[
(H_\beta f)(x)
=\frac{(-\log x)^\beta}{x}
\int_0^x f(t)(-\log t)^{-\beta}\,dt,
\qquad 0<x<1,
\]
for $\beta>-1/2$, to adjoint shifts on
Hardy and weighted Bergman spaces. 
Using Shimorin's results, they established the 
weak approximate spectral synthesis for $H_\beta$ when
$0\le\beta\le1/2$
\cite[Theorem~3.6]{FGY}. Related results on approximate spectral synthesis
have also been obtained for certain Dirichlet
type spaces \cite{Yilmaz2018,Yilmaz}.

 The radial kernels considered below have
been studied in the theory of analytic
Hilbert spaces and unilateral weighted shifts;
see \cite{Shields1974,Shimorin2001}.  
We consider radial reproducing kernels on
$\D$ of the form
\begin{equation}\label{eq:reference}
k(z,w)=\sum_{n=0}^{\infty}
\frac{(z\overline w)^n}{w_n},
\qquad w_0=1,\quad w_n>0.
\end{equation}
Our basic assumption is that the coordinate
multiplier $S_{\mathcal H(k)}$ is a hyponormal
contraction. 
 
\begin{theorem}\label{thm:reference}
Let $k$ be a radial reproducing kernel on $\D$ of the form
\eqref{eq:reference}. If the coordinate multiplier
$S_{\mathcal H(k)}$ is a hyponormal contraction, then
$\mathcal H(k)$ has the weak approximate spectral
synthesis property.
\end{theorem}

We next transfer the   result to more general reproducing kernels. For positive kernels $k$ and $K$ on $\D$, write $k\preceq K$ when $K-k$ is positive semidefinite on every finite set. Aronszajn's kernel inclusion theorem \cite{Aronszajn} shows that the comparison $c\,k\preceq K\preceq C\,k$ identifies the corresponding reproducing kernel spaces as function spaces with equivalent Hilbert norms.  

\begin{theorem}\label{thm:main}
Let $\mathcal H(K)\subset\operatorname{Hol}(\D)$
be a reproducing kernel Hilbert space.
Suppose that there exists a radial kernel $k$
of the form \eqref{eq:reference} for which
$S_{\mathcal H(k)}$ is a hyponormal contraction,
and that there are constants $0<c\le C<\infty$
such that
\begin{equation}\label{eq:comparison}
 c\,k\preceq K\preceq C\,k.
\end{equation}
Then the shift operator $S_{\mathcal H(K)}$
is bounded.
Moreover, $\mathcal H(K)$ has the weak approximate
spectral synthesis property.
\end{theorem}

Consider the kernels
\begin{equation}\label{eq:standardkernel}
k_\alpha(z,w)
=(1-z\overline w)^{-(\alpha+2)},
\qquad \alpha>-2.
\end{equation}
The space $\mathcal H(k_{-1})$ is the Hardy
space $H^2$, while $\mathcal H(k_\alpha)$
is the standard weighted Bergman space
$A_\alpha^2$ for $\alpha>-1$.

It is well known that the coordinate
multipliers on $H^2$ and $A_\alpha^2$
are subnormal contractions;
see \cite{Sarason1998,Shapiro1977}.
Since every subnormal operator is hyponormal,
and the kernels $k_\alpha$ are radial,
these spaces satisfy the hypotheses of
Theorem~\ref{thm:reference}.
Consequently, we obtain the following result.

\begin{corollary}\label{cor:standard}
The Hardy space $H^2$ and every standard
weighted Bergman space $A_\alpha^2$,
$\alpha>-1$, have the weak approximate
spectral synthesis property.
\end{corollary}

This completely resolves the weak approximate
spectral synthesis problem for the standard
weighted Bergman spaces, as discussed by
Shimorin \cite[\S5]{Shimorin2000}.

  Let $\mu$ be a finite nonzero
positive Borel measure on $\D$ invariant under rotations such that
$\mu(\{z\in\D:|z|>r\})>0$ for every $0<r<1$.
The corresponding weighted Bergman space is
\begin{equation}\label{eq:radialmeasureBergman}
 A^2(\mu)
 =\left\{f\in\operatorname{Hol}(\D):
 \|f\|_{A^2(\mu)}^2
 :=\int_\D |f(z)|^2\,d\mu(z)<\infty\right\}.
\end{equation}

Since $\mu$ is rotation-invariant, the
corresponding reproducing kernel is radial.
Moreover, the coordinate multiplier on
$A^2(\mu)$ is a subnormal contraction
\cite[p.~246]{Shapiro1977}, and hence
a hyponormal contraction.
Thus Theorem~\ref{thm:reference} yields
the following result.

\begin{corollary}\label{cor:radialmeasure}
Let $\mu$ be a finite nonzero
positive Borel measure invariant under rotations
on $\D$ such that
\[
\mu(\{z\in\D:|z|>r\})>0
\qquad(0<r<1).
\]
Then the Bergman space $A^2(\mu)$ defined
in \eqref{eq:radialmeasureBergman} has the
weak approximate spectral synthesis property.

\end{corollary}

This gives an affirmative answer to the
question raised by Shimorin
\cite[\S5]{Shimorin2000}
on extending weak approximate spectral
synthesis to general radial weighted
Bergman spaces.

As a further application, we consider the
logarithmic Hardy operators studied by
Fu, Guo and Yan \cite{FGY}. We establish
weak approximate spectral synthesis for
$H_\beta$ for all $\beta\ge0$, while the
remaining case $-1/2<\beta<0$ is not covered
by our results and, to the best of our
knowledge, remains open.

\begin{corollary}\label{cor:Hardyoperator}
For every $\beta\ge0$ and every closed subspace $N\subset L^2[0,1]$ invariant under $H_\beta$, there exist finite sets
$\Sigma_j\subset\{s\in\mathbb C:\operatorname{Re}s>-1/2\}$ such that
\begin{equation}\label{eq:Hardyapplication}
 N=\liminf_{j\to\infty}
 \operatorname{span}\{x^s(-\log x)^\beta:s\in\Sigma_j\}.
\end{equation}
The elements of each $\Sigma_j$ can be chosen distinct.
\end{corollary}

Our approach differs substantially from
that of Shimorin \cite{Shimorin2000}.
Instead of relying on extremal functions
and the positivity of weighted biharmonic
Green functions, we exploit the hyponormality
of the coordinate multiplier through the
logarithmic convexity of the squared monomial
norms.
This enables us to construct polynomial
approximants with finitely many simple zeros
whose multiplication operators are uniformly
bounded below.

The key step is to show that every weak limit
of bounded sequences from the corresponding
finite-codimensional invariant subspaces belongs
to a prescribed cyclic invariant subspace.
This yields weak approximate spectral synthesis
for the orthogonal complements of cyclic
invariant subspaces. We then extend the result
to arbitrary adjoint-invariant subspaces using
a countable family of generators and a diagonal
selection argument. Finally, comparison of reproducing kernels transfers the result to more general
analytic Hilbert spaces.

Section~\ref{sec:prelim} introduces the radial shift and the relevant weight estimates. Section~\ref{sec:division} establishes the uniform result on polynomial division, and Section~\ref{sec:referenceproof} proves the radial synthesis theorem. Kernel comparison and applications are developed in  Section~\ref{sec:applications}.

\section{Preliminaries}
\label{sec:prelim}

Let $\mathcal H(k)$ be the reproducing kernel
Hilbert space associated with the radial kernel
\eqref{eq:reference}.
The realization of such spaces as weighted
spaces of analytic functions is standard in
the theory of unilateral weighted shifts;
see \cite{Shields1974} and
\cite[Section~4]{Shimorin2001}.
The monomials form an orthogonal basis of
$\mathcal H(k)$, and
\[
\left\|\sum_{n=0}^{\infty}a_nz^n
\right\|_{\mathcal H(k)}^2
=
\sum_{n=0}^{\infty}|a_n|^2w_n.
\]
In particular,
\[
e_n(z)=\frac{z^n}{\sqrt{w_n}},
\qquad n\ge0,
\]
is an orthonormal basis of $\mathcal H(k)$,
and polynomials are dense in $\mathcal H(k)$.

The coordinate multiplier satisfies
\[
S_{\mathcal H(k)}e_n
=
\left(\frac{w_{n+1}}{w_n}\right)^{1/2}
e_{n+1}.
\]
By the classical criterion for unilateral
weighted shifts \cite{Shields1974},
$S_{\mathcal H(k)}$ is a hyponormal contraction
if and only if the ratios $w_{n+1}/w_n$
are nondecreasing and bounded above by one.
Moreover, since
\[
S_{\mathcal H(k)}^*k_\lambda
=\overline{\lambda}k_\lambda,
\qquad \lambda\in\D,
\]
we have $\|S_{\mathcal H(k)}\|\ge1$
whenever $S_{\mathcal H(k)}$ is bounded.
Consequently, if $S_{\mathcal H(k)}$ is a
hyponormal contraction, then
\[
\lim_{n\to\infty}\frac{w_{n+1}}{w_n}
=\|S_{\mathcal H(k)}\|^2=1.
\]
Thus $S_{\mathcal H(k)}$ is a hyponormal
contraction if and only if
\begin{equation}\label{eq:weightcondition}
\frac{w_{n+1}}{w_n}\nearrow1.
\end{equation}
Here $\nearrow$ denotes nondecreasing convergence.
In the proofs, we shall work with
the equivalent sequence condition
\eqref{eq:weightcondition}.

The monotonicity of $w_{n+1}/w_n$
also implies that $(w_n)$ is
logarithmically convex:
\[
w_{n+1}^2\le w_nw_{n+2},
\qquad n\ge0.
\]
Moreover, since $w_{n+1}/w_n\le1$,
we have $w_n\downarrow c\ge0$.
If $c>0$, then
\[
c\|f\|_{H^2}^2
\le
\|f\|_{\mathcal H(k)}^2
\le
\|f\|_{H^2}^2.
\]
Thus $\mathcal H(k)$ and $H^2$
coincide as spaces of analytic functions,
with equivalent Hilbert norms.
For the polynomial perturbation method,
we shall need the following estimates
in the case $c=0$.

\begin{lemma}\label{lem:weights}
Let $(w_n)$ be as in \eqref{eq:reference}
and satisfy \eqref{eq:weightcondition}.
Suppose that $w_n\to0$.
Then, for all $n,k\ge0$,
\begin{equation}\label{eq:weightfacts}
w_{n+k}\ge w_nw_k,
\qquad
\frac{w_{n+k}}{w_n}\longrightarrow1
\quad(n\to\infty,\ k\text{ fixed}).
\end{equation}
Moreover, for every fixed integer $d\ge0$,
\begin{equation}\label{eq:And}
\frac{w_{2n-d}}{w_nw_{n-d}}
\longrightarrow\infty
\qquad(n\to\infty).
\end{equation}
\end{lemma}

\begin{proof}
Since $w_0=1$ and the ratios
$w_{n+1}/w_n$ are nondecreasing,
\[
\frac{w_{n+k}}{w_n}
=
\prod_{j=0}^{k-1}
\frac{w_{n+j+1}}{w_{n+j}}
\ge
\prod_{j=0}^{k-1}
\frac{w_{j+1}}{w_j}
=w_k.
\]
For each fixed $k$, the product on the
left tends to one, proving
\eqref{eq:weightfacts}.

Fix integers $d,k\ge0$.
By logarithmic convexity, for
$n\ge d+k$,
\[
w_nw_{n-d}
\le w_kw_{2n-d-k}.
\]
Consequently,
\[
\limsup_{n\to\infty}
\frac{w_nw_{n-d}}{w_{2n-d}}
\le
w_k\lim_{n\to\infty}
\frac{w_{2n-d-k}}{w_{2n-d}}
=w_k.
\]
Since $w_k\to0$ as $k\to\infty$,
we obtain \eqref{eq:And}.
\end{proof}

\section{Proof of Theorem~\ref{thm:reference}}
\label{sec:referenceproof}

\subsection{Uniform Lower Bounds}
\label{sec:division}

\begin{proposition}\label{prop:perturbation}

Assume that $w_n\to0$.
Given a polynomial $v$ and $\varepsilon>0$,
there exists an integer
$N>\max\{0,\deg v\}$ such that,
for every integer $n\ge N$,
one can choose $\delta\in\mathbb C$ with
$|\delta|<\min\{1,\varepsilon\}$ such that
\begin{equation}\label{eq:perturbation}
 p(z)=v(z)+\frac{z^n}{\sqrt{w_n}}+\delta
\end{equation}
has $n$ distinct zeros in $\D$ and
\begin{equation}\label{eq:divisionestimate}
\norm{pq}_{\mathcal H(k)}^2
\ge\frac{1}{2w_n}
\norm{z^nq}_{\mathcal H(k)}^2
\ge\frac12\norm q_{\mathcal H(k)}^2
\qquad(q\in\mathcal H(k)).
\end{equation}
If $Z(p)$ is the zero set of $p$, then
\begin{equation}\label{eq:zeroideal}
 p\mathcal H(k)=\{h\in \mathcal H(k):h|_{Z(p)}=0\},\qquad
 (p\mathcal H(k))^\perp=\operatorname{span}\{k_\lambda:\lambda\in Z(p)\}.
\end{equation}
\end{proposition}

\begin{proof}
Put $d=\max\{0,\deg v\}$, with
$\deg 0=-\infty$, and write
\[
v(z)=\sum_{\ell=0}^{d}v_\ell z^\ell.
\]
Let
\[
V=1+\sum_{\ell=0}^{d}|v_\ell|.
\]

Since $S_{\mathcal H(k)}$ is a contraction,
the multiplication operator $M_{v+\delta}$
on $\mathcal H(k)$ satisfies, for $|\delta|<1$, 
\[
\norm{M_{v+\delta}}\le V,
\qquad
\norm{v+\delta}_\infty\le V.
\] 
For $n>d$, set
\[
A_n(d)^2
:=
\frac{w_{2n-d}}{w_nw_{n-d}}.
\]
By Lemma~\ref{lem:weights}, every sufficiently
large integer $n>d$ satisfies
\begin{equation}\label{eq:nchoice}
A_n(d)\ge4V,
\qquad
w_n^{-1/2}>2V.
\end{equation}

For $q(z)=\sum_{m\ge0}q_mz^m$, let
\[
(R_{n-d}q)(z)
=\sum_{m=n-d}^{\infty}q_mz^m.
\]
Since $\deg(v+\delta)\le d$, the polynomial
$(v+\delta)(q-R_{n-d}q)$ has degree less
than $n$, whereas $e_nq$ contains only
terms of degree at least $n$.
By orthogonality,
\begin{equation}\label{eq:cross}
\ip{e_nq}{(v+\delta)q}_{\mathcal H(k)}
=
\ip{e_nq}{(v+\delta)R_{n-d}q}_{\mathcal H(k)}.
\end{equation}

For fixed $n$, the ratio
\[
\frac{w_{n+m}}{w_m}
=
\prod_{j=0}^{n-1}
\frac{w_{m+j+1}}{w_{m+j}}
\]
is nondecreasing in $m$.
Thus, for $m\ge n-d$,
\[
\frac{w_{n+m}}{w_n}
\ge
\frac{w_{2n-d}}{w_nw_{n-d}}w_m
=A_n(d)^2w_m.
\]
Consequently,
\begin{align}
\norm{e_nq}_{\mathcal H(k)}^2
&=
\sum_{m\ge0}|q_m|^2
\frac{w_{n+m}}{w_n}
\notag\\
&\ge
A_n(d)^2
\sum_{m\ge n-d}|q_m|^2w_m
\notag\\
&=
A_n(d)^2
\norm{R_{n-d}q}_{\mathcal H(k)}^2.
\label{eq:tailcontrol}
\end{align}

Combining \eqref{eq:cross} and
\eqref{eq:tailcontrol}, we obtain
\[
\begin{aligned}
\left|
\ip{e_nq}{(v+\delta)q}_{\mathcal H(k)}
\right|
&\le
\norm{e_nq}_{\mathcal H(k)}
\norm{(v+\delta)R_{n-d}q}_{\mathcal H(k)}\\
&\le
V\norm{e_nq}_{\mathcal H(k)}
\norm{R_{n-d}q}_{\mathcal H(k)}\\
&\le
\frac{V}{A_n(d)}
\norm{e_nq}_{\mathcal H(k)}^2.
\end{aligned}
\]
It follows from \eqref{eq:nchoice} that
\[
\begin{aligned}
\norm{pq}_{\mathcal H(k)}^2
&=
\norm{e_nq+(v+\delta)q}_{\mathcal H(k)}^2\\
&\ge
\left(1-\frac{2V}{A_n(d)}\right)
\norm{e_nq}_{\mathcal H(k)}^2\\
&\ge
\frac12\norm{e_nq}_{\mathcal H(k)}^2\\
&=
\frac{1}{2w_n}
\norm{z^nq}_{\mathcal H(k)}^2.
\end{aligned}
\]
Moreover, by \eqref{eq:weightfacts},
\[
\frac1{w_n}\norm{z^nq}_{\mathcal H(k)}^2
=
\sum_{m\ge0}|q_m|^2
\frac{w_{n+m}}{w_n}
\ge
\norm q_{\mathcal H(k)}^2.
\]
This proves \eqref{eq:divisionestimate}.

We next choose $\delta$ so that $p$
has only simple zeros. A multiple zero
$\zeta$ of $p=e_n+v+\delta$ must satisfy
\[
v'(\zeta)
+\frac{n\zeta^{n-1}}{\sqrt{w_n}}=0,
\qquad
\delta=-v(\zeta)
-\frac{\zeta^n}{\sqrt{w_n}}.
\]
Since $n>d$, the derivative of $v+e_n$
is a nonzero polynomial and hence
has only finitely many zeros.
Thus there are only finitely many
values of $\delta$ for which $p$
has a multiple zero.
We may therefore choose
$|\delta|<\min\{1,\varepsilon\}$
outside this finite set.

For $|z|=1$, \eqref{eq:nchoice} gives
\[
|e_n(z)|=w_n^{-1/2}
>2V>|v(z)+\delta|.
\]
By Rouch\'e's theorem, $p$ and $e_n$
have the same number of zeros in $\D$,
counted with multiplicity.
Since $p$ has degree $n$ and no
multiple zeros, all its $n$ zeros
are distinct and belong to $\D$.

Finally, \eqref{eq:divisionestimate}
implies that $p\mathcal H(k)$ is closed.
For every polynomial $h$, Euclidean
division by $p$ gives
\[
h=pa+r,\qquad \deg r<n,
\]
where $a$ and $r$ are polynomials.
Since polynomials are dense in
$\mathcal H(k)$, their residue classes
are dense in the quotient space
$\mathcal H(k)/p\mathcal H(k)$.
The division formula shows that these
classes lie in the span of the residue classes of
$1,z,\ldots,z^{n-1}$.
This span has finite dimension and is closed, and hence
\[
\operatorname{codim}p\mathcal H(k)\le n.
\]

On the other hand, let
\[
Z(p)=\{\lambda_1,\ldots,\lambda_n\}
\]
and define
\[
\mathcal V_p=
\{h\in\mathcal H(k):
h(\lambda_j)=0,\ 1\le j\le n\}.
\]
The evaluation functionals at the
distinct points $\lambda_1,\ldots,\lambda_n$
are linearly independent. Therefore,
\[
\operatorname{codim}\mathcal V_p=n.
\]
Since
\[
p\mathcal H(k)\subseteq\mathcal V_p
\]
and $\operatorname{codim}p\mathcal H(k)\le n$,
the two closed subspaces must coincide.
Thus
\[
p\mathcal H(k)=\mathcal V_p.
\]
Taking orthogonal complements and using
the reproducing property, we obtain
\[
(p\mathcal H(k))^\perp
=
\operatorname{span}
\{k_\lambda:\lambda\in Z(p)\}.
\]
This proves \eqref{eq:zeroideal}.
\end{proof}

\subsection{From Cyclic Subspaces to General Invariant Subspaces}

The next proposition uses only contractivity of the coordinate multiplier and the unitary power dilation theorem \cite[Chapter~I]{SzNagyFoias}.

For $f\in\mathcal H(k)$, let
\[
[f]=\overline{\operatorname{span}}
\{z^nf:n\ge0\}
\]
denote the shift-invariant subspace
generated by $f$.

\begin{proposition}\label{prop:total}
Let $\mathcal H(k)$ be an  analytic reproducing kernel Hilbert space on $\D$ for which the shift operator $S_{\mathcal H(k)}$ is a contraction. 
Every nonzero closed
subspace $M$ invariant under $S_{\mathcal H(k)}$
contains a sequence $(f_i)$
whose linear span is dense in $M$,
such that multiplication by each $f_i$
defines a bounded operator
from $H^2$ into $\mathcal H(k)$. 
\end{proposition}

\begin{proof}
Let $P_M$ denote the orthogonal projection
onto $M$. For every $g\in M$ and
$\lambda\in\D$, we have
\[
\begin{aligned}
\ip{g}{(S_{\mathcal H(k)}|_M)^*P_Mk_\lambda}
 =\ip{S_{\mathcal H(k)}g}{k_\lambda} =\lambda g(\lambda) =\ip{g}{\overline{\lambda}P_Mk_\lambda}.
\end{aligned}
\]
Hence
\[
(S_{\mathcal H(k)}|_M)^*P_Mk_\lambda
=\overline{\lambda}P_Mk_\lambda.
\]

Let $\mathbb T=\partial\D$, and let $m$
denote normalized arc length measure
on $\mathbb T$. 
Fix $\lambda\in\D$ such that
$P_Mk_\lambda\ne0$.
Since $S_{\mathcal H(k)}|_M$ is a
contraction, let $U$ be its unitary
power dilation. Let $\mu_\lambda$
be the spectral measure of $U$
associated with $P_Mk_\lambda$.

For every $n\ge0$,
\[
\begin{aligned}
\int_{\mathbb T}\zeta^n\,d\mu_\lambda(\zeta)
&=\ip{U^nP_Mk_\lambda}{P_Mk_\lambda}\\
&=\ip{(S_{\mathcal H(k)}|_M)^n
P_Mk_\lambda}{P_Mk_\lambda}\\
&=\lambda^n\|P_Mk_\lambda\|^2.
\end{aligned}
\]
By uniqueness of Fourier coefficients
of finite measures,
\[
d\mu_\lambda(\zeta)
=\|P_Mk_\lambda\|^2
\frac{1-|\lambda|^2}
{|\zeta-\lambda|^2}\,dm(\zeta).
\]

Consequently, for every analytic
polynomial $q$,
\[
\begin{aligned}
\|qP_Mk_\lambda\|_{\mathcal H(k)}^2
&=\|q(S_{\mathcal H(k)}|_M)
P_Mk_\lambda\|_{\mathcal H(k)}^2\\
&\le \|q(U)P_Mk_\lambda\|^2\\
&=\int_{\mathbb T}|q(\zeta)|^2
\,d\mu_\lambda(\zeta)\\
&\le
\frac{1+|\lambda|}{1-|\lambda|}
\|P_Mk_\lambda\|^2
\|q\|_{H^2}^2.
\end{aligned}
\]
By density of polynomials in $H^2$,
the map $q\mapsto qP_Mk_\lambda$
extends to a bounded operator from
$H^2$ into $\mathcal H(k)$.
Continuity of point evaluations
shows that this extension agrees
with pointwise multiplication.

Finally, choose a countable dense
subset $\{\lambda_j\}_{j\ge1}$ of
$\D$. If $g\in M$ is orthogonal to
every $P_Mk_{\lambda_j}$, then
\[
g(\lambda_j)
=\ip{g}{P_Mk_{\lambda_j}}=0
\qquad(j\ge1).
\]
By the identity theorem, $g=0$.
Thus
\[
M=\overline{\operatorname{span}}
\{P_Mk_{\lambda_j}:j\ge1\}.
\]
Discarding the zero vectors gives
the required sequence.
\end{proof}

Let $(f_i)$ be the sequence obtained above.
Since $M$ is invariant under the shift and the linear
span of $(f_i)$ is dense in $M$, we have
\[
M=\overline{\operatorname{span}}
\bigcup_{i\ge1}[f_i],
\qquad
M^\perp=\bigcap_{i\ge1}[f_i]^\perp,
\] 
Moreover, each multiplication operator
$M_{f_i}:H^2\to\mathcal H(k)$ is bounded. 
It therefore suffices to establish the desired
approximation for each $[f_i]^\perp$.

\begin{proposition}\label{prop:cyclic}
Assume that \eqref{eq:weightcondition}
holds and $w_n\to0$.
Let $f\in\mathcal H(k)$ be such that
multiplication by $f$ defines a bounded
operator $M_f:H^2\to\mathcal H(k)$.
Then there exist finite sets $A_j\subset\D$
such that
\[
[f]^\perp
=
\liminf_{j\to\infty}
\operatorname{span}
\{k_\lambda:\lambda\in A_j\}.
\]
\end{proposition}

\begin{proof}
Choose polynomials $v_j$ converging to $f$
in $\mathcal H(k)$.
By Proposition~\ref{prop:perturbation},
we may recursively choose strictly increasing
integers $n_j$ and complex numbers $\delta_j$
with $|\delta_j|<1/j$ such that
\begin{equation}\label{eq:pjsequence}
p_j=v_j+\frac{z^{n_j}}{\sqrt{w_{n_j}}}
+\delta_j,
\qquad
M_j=p_j\mathcal H(k).
\end{equation}
Each $M_j$ is closed and has finite codimension
by Proposition~\ref{prop:perturbation}.
Moreover, $M_j$ is invariant under the shift and
$p_j\in M_j$.

Since the normalized monomials form an
orthonormal sequence, $(e_{n_j})$ converges
weakly to zero.
Together with the norm convergence of $(v_j)$
to $f$ and $\delta_j\to0$, this shows that
$(p_j)$ converges weakly to $f$ in $\mathcal H(k)$.

We first show that every weak limit of a
sequence $g_j\in M_j$, taken along any
subsequence, belongs to $[f]$.
Consider an arbitrary subsequence, retaining
the index $j$.
Suppose that
\[
g_j=p_jq_j\in M_j,\qquad
\norm{g_j}_{\mathcal H(k)}\le B,
\]
and that $(g_j)$ converges weakly to
some $g\in\mathcal H(k)$.

By the second inequality in
\eqref{eq:divisionestimate},
\[
\norm{q_j}_{\mathcal H(k)}^2
\le
2\norm{p_jq_j}_{\mathcal H(k)}^2
\le 2B^2.
\]
Thus $(q_j)$ is bounded in $\mathcal H(k)$.
Passing to a further subsequence, we may
assume that $(q_j)$ converges weakly to
some $q\in\mathcal H(k)$.

On the other hand, the first inequality
in \eqref{eq:divisionestimate} gives
\[
\sum_{k\ge0}|(q_j)_k|^2
\frac{w_{n_j+k}}{w_{n_j}}
=
\frac{\norm{z^{n_j}q_j}_{\mathcal H(k)}^2}
{w_{n_j}}
\le 2B^2.
\]
Consequently, for every fixed integer $K\ge0$,
\[
\sum_{k=0}^K |(q_j)_k|^2
\frac{w_{n_j+k}}{w_{n_j}}
\le 2B^2.
\]
Since the coefficient functionals are continuous,
$(q_j)_k$ converges to $q_k$ for each fixed $k$.
Moreover, by \eqref{eq:weightfacts},
\[
\frac{w_{n_j+k}}{w_{n_j}}
\longrightarrow1
\qquad(j\to\infty)
\]
for every fixed $k$.
Letting $j\to\infty$ gives
\[
\sum_{k=0}^K|q_k|^2\le2B^2.
\]
Since $K$ is arbitrary, we conclude that
\begin{equation}\label{eq:hardylimit}
q\in H^2,
\qquad
\norm{q}_{H^2}^2\le2B^2.
\end{equation}

Since $(g_j)$, $(p_j)$, and $(q_j)$
converge weakly to $g$, $f$, and $q$,
respectively, continuity of point evaluations
gives
\[
g(z)=f(z)q(z),
\qquad z\in\D.
\]
Since $M_f:H^2\to\mathcal H(k)$ is bounded,
polynomial approximation in $H^2$ implies
that $fq\in[f]$.
Thus every weak limit obtained in this way
belongs to $[f]$.

For each $j$, let $P_j$ denote the orthogonal
projection from $\mathcal H(k)$ onto $M_j$.
Since
\[
\operatorname{dist}_{\mathcal H(k)}
(x,M_j^\perp)=\norm{P_jx}_{\mathcal H(k)},
\]
we first prove that
\[
[f]^\perp\subseteq\liminf_jM_j^\perp.
\]
Suppose, to the contrary, that
$x\perp[f]$ but $\norm{P_jx}_{\mathcal H(k)}$
does not tend to zero.
Then there exist $\varepsilon>0$ and
a subsequence $(j_\ell)$ such that
\[
\norm{P_{j_\ell}x}_{\mathcal H(k)}
\ge\varepsilon
\qquad(\ell\ge1).
\]
Since $(P_{j_\ell}x)$ is bounded, we may
pass to a further subsequence converging
weakly to some $h\in\mathcal H(k)$.
By the property of weak limits established
above, $h\in[f]$.
Since $x\perp[f]$, we obtain
\[
\norm{P_{j_\ell}x}_{\mathcal H(k)}^2
=
\ip{x}{P_{j_\ell}x}_{\mathcal H(k)}
\longrightarrow
\ip{x}{h}_{\mathcal H(k)}=0,
\]
a contradiction.
Hence
\[
[f]^\perp\subseteq\liminf_jM_j^\perp.
\]

Conversely, let $x\in\liminf_jM_j^\perp$.
By the distance identity above,
$\norm{P_jx}_{\mathcal H(k)}\to0$.
Since $(p_j)$ converges weakly to $f$
and $S_{\mathcal H(k)}$ is bounded,
for every fixed $m\ge0$,
\[
\sup_j
\norm{S_{\mathcal H(k)}^mp_j}_{\mathcal H(k)}
<\infty.
\]
Moreover, $S_{\mathcal H(k)}^mp_j\in M_j$
because $M_j$ is invariant under the shift.
Therefore,
\[
\begin{aligned}
 \left|
\ip{x}{S_{\mathcal H(k)}^mp_j}_{\mathcal H(k)}
\right| 
& =
\left|
\ip{P_jx}{S_{\mathcal H(k)}^mp_j}_{\mathcal H(k)}
\right|\\
&\quad\le
\norm{P_jx}_{\mathcal H(k)}
\sup_j\norm{S_{\mathcal H(k)}^mp_j}_{\mathcal H(k)}
\longrightarrow0.
\end{aligned}
\]
Since $S_{\mathcal H(k)}^m$ is bounded,
the sequence $(S_{\mathcal H(k)}^mp_j)$
converges weakly to $S_{\mathcal H(k)}^mf$.
Thus
\[
\ip{x}{S_{\mathcal H(k)}^mf}_{\mathcal H(k)}
=0
\qquad(m\ge0).
\]
Hence $x\perp[f]$, proving
\[
\liminf_jM_j^\perp\subseteq[f]^\perp.
\]
Combining the two inclusions, we obtain
\[
[f]^\perp=\liminf_{j\to\infty}M_j^\perp.
\]

Finally, Proposition~\ref{prop:perturbation}
gives
\[
M_j^\perp
=
\operatorname{span}
\{k_\lambda:\lambda\in Z(p_j)\},
\]
where $Z(p_j)$ consists of finitely many
distinct points of $\D$.
Taking $A_j=Z(p_j)$ completes the proof.
\end{proof}

To pass from orthogonal complements of cyclic shift-invariant subspaces to arbitrary
invariant subspaces, we use the following
lemma on diagonal selection due to Shimorin
\cite[Lemma~1]{Shimorin2000}.

\begin{lemma}\label{lem:intersection}
Let $Z$ be a separable Banach space, and let
$(E_{n,k})_{k\ge1}$ be a sequence of subspaces
of $Z$ for each $n\ge1$. Suppose that
\[
N_n=\liminf_{k\to\infty}E_{n,k},
\qquad n\ge1.
\]
Then there exist sequences of positive integers
$(n_j)$ and $(k_j)$ such that
\[
\bigcap_{n=1}^{\infty}N_n
=
\liminf_{j\to\infty}E_{n_j,k_j}.
\]
\end{lemma}

\begin{proof}[Proof of Theorem~\ref{thm:reference}]
By the monotonicity of $w_{n+1}/w_n$
and the fact that these ratios converge
to one, we have $w_n\downarrow c\ge0$. 
Let $N\in\operatorname{Lat}(S_{\mathcal H(k)}^*)$.

We first consider the trivial cases.
If $N=\{0\}$, take $E_j=\{0\}$ for every $j$.
If $N=\mathcal H(k)$, choose a sequence of distinct points
$(\lambda_j)$ in $\D$ having an accumulation
point in $\D$, and set
\[
E_j=\operatorname{span}
\{k_{\lambda_1},\ldots,k_{\lambda_j}\}.
\]
Each $E_j$ is finite dimensional and
invariant under $S_{\mathcal H(k)}^*$.
If $f\in \mathcal H(k)$ is orthogonal to every $E_j$, then
\[
f(\lambda_j)
=\langle f,k_{\lambda_j}\rangle_{\mathcal H(k)}=0
\qquad(j\ge1).
\]
By the identity theorem, $f=0$. Hence
\[
\overline{\bigcup_{j\ge1}E_j}=\mathcal H(k).
\]
Since $(E_j)$ is increasing, it follows that
\[
\mathcal H(k)=\liminf_{j\to\infty}E_j.
\]

We may therefore assume that
$\{0\}\subsetneq N\subsetneq \mathcal H(k)$.

Suppose first that $c=0$.
Since $N$ is invariant under $S_{\mathcal H(k)}^*$, its orthogonal
complement $N^\perp$ is a nonzero closed
subspace invariant under $S_{\mathcal H(k)}$.

By Proposition~\ref{prop:total}, there exists
a countable family $(f_i)\subset N^\perp$
whose linear span is dense in $N^\perp$
and such that, for each $i$, multiplication
by $f_i$ defines a bounded operator
from $H^2$ into $\mathcal H(k)$.

Since $N^\perp$ is invariant under $S_{\mathcal H(k)}$,
\[
N^\perp
=\overline{\operatorname{span}}\bigcup_i[f_i].
\]
Taking orthogonal complements gives
\[
N=\bigcap_i[f_i]^\perp.
\]

By Proposition~\ref{prop:cyclic}, each
$[f_i]^\perp$ is the lower limit of a sequence
of subspaces of finite dimension invariant
under $S_{\mathcal H(k)}^*$ and spanned by reproducing kernels
at distinct points of $\D$. Then 
Lemma~\ref{lem:intersection} yields the
desired representation of $N$ as a lower limit
of finite-dimensional subspaces spanned by reproducing kernels.

Suppose now that $c>0$.
Since
\[
c\le w_n\le1,
\]
we have
\[
c\norm f_{H^2}^2
\le\norm f_{\mathcal H(k)}^2
\le\norm f_{H^2}^2.
\]
Thus $\mathcal H(k)$ and $H^2$ coincide as spaces of
holomorphic functions, with equivalent
Hilbert norms.

Let $J:H^2\to \mathcal H(k)$ be the identity map.
Then $J$ and $J^*$ are bounded and invertible,
and
\[
S_{\mathcal H(k)}J=JS_{H^2},
\qquad
J^*S_{\mathcal H(k)}^*=S_{H^2}^*J^*.
\]

Set $N_0=J^*N$. Since $J^*$ is bounded
and invertible, $N_0$ is closed. Moreover,
\[
S_{H^2}^*N_0
=S_{H^2}^*J^*N
=J^*S_{\mathcal H(k)}^*N
\subseteq J^*N=N_0.
\]
Thus $N_0$ is a closed
subspace invariant under $S_{H^2}^*$.

By the classical approximation
result for the Hardy space \cite[\S4.3]{DSS}
(see also \cite[p.~4]{Shimorin2000}),
together with approximation by finite
Blaschke products having simple zeros,
there exist finite sets $A_j\subset\D$
of distinct points such that
\[
N_0=\liminf_{j\to\infty}F_j,
\qquad
F_j=\operatorname{span}
\{k_\lambda^{H^2}:\lambda\in A_j\}.
\]

For every $f\in H^2$ and $\lambda\in\D$,
the reproducing property gives
\[ 
\langle f,J^*k_\lambda\rangle_{H^2}
 =\langle Jf,k_\lambda\rangle_{\mathcal H(k)} =f(\lambda) =\langle f,k_\lambda^{H^2}\rangle_{H^2}. 
\]
Hence
\[
J^*k_\lambda=k_\lambda^{H^2}.
\]
Set
\[
E_j=(J^*)^{-1}F_j
=\operatorname{span}
\{k_\lambda:\lambda\in A_j\}.
\]

For every $x\in \mathcal H(k)$, bounded invertibility
of $J^*$ yields
\[
\begin{aligned}
\frac{1}{\norm{J^*}}
\operatorname{dist}_{H^2}(J^*x,F_j)
&\le\operatorname{dist}_{\mathcal H(k)}(x,E_j)\\
&\le\norm{(J^*)^{-1}}
\operatorname{dist}_{H^2}(J^*x,F_j).
\end{aligned}
\]

The preceding inequalities imply that
\[
\operatorname{dist}_{\mathcal H(k)}(x,E_j)\to0
\quad\text{if and only if}\quad
\operatorname{dist}_{H^2}(J^*x,F_j)\to0.
\]
Since $N_0=\liminf_{j\to\infty}F_j$ and
$N_0=J^*N$, we conclude that
\[
N=\liminf_{j\to\infty}E_j.
\]
This completes the proof.
\end{proof}

\section{Kernel comparison  and applications}\label{sec:applications}

\subsection{Kernel Comparison}

\begin{proof}[Proof of Theorem~\ref{thm:main}]

Let $J:\mathcal H(k)\to\mathcal H(K)$
be the identity map. 
By Aronszajn's reproducing kernel inclusion theorem
\cite[Part~I, \S7]{Aronszajn}, the kernel comparison
\eqref{eq:comparison} implies that $\mathcal H(K)$ and
$\mathcal H(k)$ coincide as spaces of holomorphic functions,
with
\begin{equation}\label{eq:normcomparison}
 C^{-1/2}\norm{f}_{\mathcal H(k)}
 \le\norm{f}_{\mathcal H(K)}
 \le c^{-1/2}\norm{f}_{\mathcal H(k)}.
\end{equation}

Consequently, $J$ and $J^*$ are bounded
and invertible, with
$\norm{J}\le c^{-1/2}$ and
$\norm{J^{-1}}\le C^{1/2}$. 
Since $S_{\mathcal H(k)}$ is a contraction,
coordinate multiplication defines a bounded
operator on $\mathcal H(K)$ satisfying
\[
S_{\mathcal H(K)}
=JS_{\mathcal H(k)}J^{-1}.
\]
In fact, for every $m\ge0$,
\[
\norm{S_{\mathcal H(K)}^m}_{\mathcal B(\mathcal H(K))}
\le\norm{J}\norm{J^{-1}}
\le\sqrt{C/c}.
\]
Taking adjoints yields 
\begin{equation}\label{eq:Jbounds}
 S_{\mathcal H(K)}J=JS_{\mathcal H(k)},\qquad J^*S_{\mathcal H(K)}^*=S_{\mathcal H(k)}^*J^*.
\end{equation}
Polynomials are dense in $\mathcal H(K)$ because
they are dense in $\mathcal H(k)$ and the two
Hilbert norms are equivalent.

Fix $N\in\operatorname{Lat}(S_{\mathcal H(K)}^*)$ and set $N_0=J^*N$.
Since $J^*$ is bounded and invertible, $N_0$ is closed.
Furthermore,
\[
 S_{\mathcal H(k)}^*N_0=S_{\mathcal H(k)}^*J^*N=J^*S_{\mathcal H(K)}^*N\subseteq J^*N=N_0.
\]
By Theorem~\ref{thm:reference}, there are finite sets
$A_j\subset\D$ of distinct points such that
\[
 N_0=\liminf_jF_j,\qquad
 F_j=\operatorname{span}\{k_\lambda:\lambda\in A_j\}.
\]

For every $f\in\mathcal H(k)$ and $\lambda\in\D$,
the reproducing property gives
\[
 \begin{aligned}
 \ip{f}{J^*K_\lambda}_{\mathcal H(k)}
 &=\ip{Jf}{K_\lambda}_{\mathcal H(K)}
 =(Jf)(\lambda)\\
 &=f(\lambda)=\ip{f}{k_\lambda}_{\mathcal H(k)}.
 \end{aligned}
\]
Therefore
\begin{equation}\label{eq:kernelidentity}
 J^*K_\lambda=k_\lambda.
\end{equation}
Define $E_j=(J^*)^{-1}F_j$. Then
\[
 E_j=\operatorname{span}\{K_\lambda:\lambda\in A_j\},
\]
so each $E_j$ has finite dimension, is invariant under $S_{\mathcal H(K)}^*$,
and is spanned by kernels at distinct interior points.

For every $x\in\mathcal H(K)$,
\[
 \operatorname{dist}_{\mathcal H(k)}(J^*x,F_j)
 =\inf_{e\in E_j}\norm{J^*(x-e)}_{\mathcal H(k)}.
\]
The bounds for $J$ and $J^{-1}$ imply
\begin{equation}\label{eq:targetdistance}
 \sqrt c\,\operatorname{dist}_{\mathcal H(k)}(J^*x,F_j)
 \le\operatorname{dist}_{\mathcal H(K)}(x,E_j)
 \le\sqrt C\,\operatorname{dist}_{\mathcal H(k)}(J^*x,F_j).
\end{equation}
Thus the two distances tend to zero simultaneously.
Since $N_0=\liminf_jF_j=J^*N$, it follows that
$N=\liminf_jE_j$, as asserted.
\end{proof}

\begin{corollary}\label{cor:similarity}
Let $\mathcal H(K)$ satisfy the hypotheses of
Theorem~\ref{thm:main}, and let $
L:\mathcal H(K)\longrightarrow\mathscr H$ 
be a bounded invertible operator between
Hilbert spaces. Set
\[
T=LS_{\mathcal H(K)}^*L^{-1}.
\]
Then every closed subspace $N\subset\mathscr H$
invariant under $T$ has a representation
\[
N=\liminf_{j\to\infty}
\operatorname{span}\{LK_\lambda:\lambda\in A_j\},
\]
where each $A_j$ is a finite subset of $\D$
consisting of distinct points and
$K_\lambda=K(\cdot,\lambda)$.

\end{corollary}

\begin{proof} 
Let $N\in\operatorname{Lat}(T)$ and put
$N_0=L^{-1}N$. Then
$N_0\in\operatorname{Lat}(S_{\mathcal H(K)}^*)$. 
By Theorem~\ref{thm:main}, there exist finite
kernel spans $F_j$ such that
$N_0=\liminf_jF_j$.

Since $L$ is bounded and invertible,
\[
\frac{\operatorname{dist}(x,F_j)}{\|L^{-1}\|}
\le \operatorname{dist}(Lx,LF_j)
\le \|L\|\operatorname{dist}(x,F_j).
\]
Consequently,
\[
N=\liminf_{j\to\infty}LF_j.
\]
Finally,
\[
T(LK_\lambda)=\overline\lambda LK_\lambda,
\]
so the spaces $LF_j$ have finite dimension,
are invariant under $T$, and are spanned by eigenvectors.
\end{proof}

\medskip
The preceding results establish weak approximate
spectral synthesis for adjoints of coordinate
multipliers. It is natural to ask whether this
property extends to more general analytic
multiplication operators.

\medskip
\noindent\textbf{Question.}
Let $\mathcal H(K)$ satisfy the hypotheses of
Theorem~\ref{thm:main}, and let
$\varphi\in\operatorname{Mult}(\mathcal H(K))$.
For which analytic multipliers $\varphi$ does
$M_\varphi^*$ possess weak approximate
spectral synthesis? 

\subsection{Radial Bergman Spaces}

We first consider the standard weighted
Bergman spaces and then extend the result
to  radial weights.

\begin{proof}[Proof of Corollary~\ref{cor:standard}]
The reproducing kernels of $H^2$ and $A_\alpha^2$
are $k_{-1}$ and $k_\alpha$, respectively.
Both kernels are radial.
Moreover, the corresponding coordinate
multipliers are subnormal contractions
\cite[p.~20]{Sarason1998}
\cite[p.~246]{Shapiro1977},
and hence hyponormal contractions.
The conclusion follows directly from
Theorem~\ref{thm:reference}.
\end{proof}

The standard weighted Bergman spaces are special
cases of the more general radial weighted Bergman
spaces studied extensively in the literature
\cite{HKZ,Shimorin2000,Shimorin2003}.

\begin{proof}[Proof of Corollary~\ref{cor:radialmeasure}]
By normalizing the measure, we may assume
that $\mu(\D)=1$.
Set
\[
w_n=\int_\D |z|^{2n}\,d\mu(z),
\qquad n\ge0.
\]
Then $w_0=1$ and $w_n>0$.
By rotational invariance and Parseval's identity,
\[
\left\|\sum_{n=0}^\infty a_nz^n
\right\|_{A^2(\mu)}^2
=
\sum_{n=0}^\infty |a_n|^2w_n.
\]
For every $0<r<1$, we have
\[
r^{2n}\mu(\{z\in\D:|z|>r\})
\le w_n\le1.
\]
The support condition on $\mu$ therefore
implies that
\[
\lim_{n\to\infty}w_n^{1/n}=1.
\]
Consequently, $A^2(\mu)$ is a reproducing
kernel Hilbert space with radial kernel
\[
K_\mu(z,\lambda)
=
\sum_{n=0}^\infty
\frac{(z\overline\lambda)^n}{w_n}.
\]

The coordinate multiplier on $A^2(\mu)$
is the restriction of the normal contraction
$M_z$ on $L^2(\mu)$ to the invariant
subspace $A^2(\mu)$.
Hence it is a subnormal contraction
\cite[p.~246]{Shapiro1977},
and therefore a hyponormal contraction.

The conclusion follows from
Theorem~\ref{thm:reference}.
\end{proof}

\begin{remark}\label{rem:shimorinradial}
The assumptions on $\mu$ in
Corollary~\ref{cor:radialmeasure} are precisely
those used by Shimorin
\cite[\S4, condition~(18)]{Shimorin2000},  
up to normalization.

For this class of radial measures, Shimorin
established a Carath\'eodory--Schur
approximation theorem for subextremal functions
\cite[Theorem~5]{Shimorin2000}.
In the concluding discussion
\cite[\S5, p.~36]{Shimorin2000}, he pointed
out that extending weak approximate spectral
synthesis to arbitrary radial weighted Bergman
spaces required a different approach, since
the associated weighted biharmonic Green
function need not be positive. 
\end{remark}

\subsection{Logarithmic Hardy Operators}

We next apply the radial synthesis theorem
to the logarithmic Hardy operators through
the coanalytic model of Fu, Guo and Yan
\cite[Theorem~3.1]{FGY}.

For $\beta>-1/2$, put $\alpha=2\beta-1$. 
The cited model provides a unitary operator
\[
\Lambda_\beta:
L^2[0,1]\longrightarrow\mathcal H(k_\alpha)
\]
satisfying
\begin{equation}\label{eq:intertwineHardy}
\Lambda_\beta H_\beta\Lambda_\beta^{-1}
=
I-S_{\mathcal H(k_\alpha)}^*.
\end{equation}
Moreover, for $\operatorname{Re}s>-1/2$, 
\begin{equation}\label{eq:eigenfunctions}
\Lambda_\beta
\bigl(x^s(-\log x)^\beta\bigr)
=
\frac{\Gamma(\alpha+2)^{1/2}}
{(1+s)^{\alpha+2}}
k^\alpha_{\overline s/(1+\overline s)}.
\end{equation}

\begin{proof}[Proof of Corollary~\ref{cor:Hardyoperator}]
For $\beta\ge0$, we have
$\alpha=2\beta-1\ge-1$. 
Thus Corollary~\ref{cor:standard} applies to
$\mathcal H(k_\alpha)$. Moreover,
\[
\operatorname{Lat}(I-S_{\mathcal H(k_\alpha)}^*)
=
\operatorname{Lat}(S_{\mathcal H(k_\alpha)}^*).
\]

Let $N\in\operatorname{Lat}(H_\beta)$.
By \eqref{eq:intertwineHardy},
$\Lambda_\beta N$ is a closed
subspace of $\mathcal H(k_\alpha)$ invariant
under $S_{\mathcal H(k_\alpha)}^*$.
Hence Corollary~\ref{cor:standard} yields
finite sets $A_j\subset\D$ such that
\[
\Lambda_\beta N
=
\liminf_{j\to\infty}E_j,
\qquad
E_j=\operatorname{span}
\{k_\lambda^\alpha:\lambda\in A_j\}.
\]

For each $\lambda\in A_j$, set
\[
s=\frac{\overline\lambda}{1-\overline\lambda}.
\]
Since $|\lambda|<1$, we have
\[
2\operatorname{Re}s+1
=
\frac{1-|\lambda|^2}{|1-\lambda|^2}>0.
\]
Thus $\operatorname{Re}s>-1/2$, and
\eqref{eq:eigenfunctions} gives
\[
\Lambda_\beta f_s
=
\frac{\Gamma(\alpha+2)^{1/2}}
{(1+s)^{\alpha+2}}
k_\lambda^\alpha,
\qquad
f_s(x)=x^s(-\log x)^\beta.
\]
The scalar factor is nonzero.

Define
\[
\Sigma_j
=
\left\{
\frac{\overline\lambda}{1-\overline\lambda}
:\lambda\in A_j
\right\}
\]
and
\[
F_j
=
\operatorname{span}\{f_s:s\in\Sigma_j\}.
\]
Since the map
$\lambda\mapsto\overline\lambda/(1-\overline\lambda)$
is injective, the elements of each
$\Sigma_j$ are distinct.
Moreover, \eqref{eq:eigenfunctions} implies
\[
\Lambda_\beta F_j=E_j.
\]

Since $\Lambda_\beta$ is unitary, for every
$f\in L^2[0,1]$ we have
\[
\operatorname{dist}_{L^2}(f,F_j)
=
\operatorname{dist}_{\mathcal H(k_\alpha)}
(\Lambda_\beta f,E_j).
\]
Consequently,
\[
\begin{aligned}
f\in N
&\iff \Lambda_\beta f\in\Lambda_\beta N\\
&\iff
\operatorname{dist}_{\mathcal H(k_\alpha)}
(\Lambda_\beta f,E_j)\to0\\
&\iff
\operatorname{dist}_{L^2}(f,F_j)\to0.
\end{aligned}
\]
Therefore,
\[
N=\liminf_{j\to\infty}F_j
=
\liminf_{j\to\infty}
\operatorname{span}
\{x^s(-\log x)^\beta:s\in\Sigma_j\},
\]
which proves \eqref{eq:Hardyapplication}.
\end{proof}

\subsection*{Author contributions}
Zhaopeng Lin, Yufeng Lu, and Chao Zu contributed equally to the conception,
development, and writing of this manuscript. All authors read and approved
the final manuscript.

\subsection*{Conflict of interest}
The authors declare no conflicts of interest.

\subsection*{Data availability statement}
No data, models, or code were generated or used for the research described in the article.

\subsection*{Funding}

Y. Lu was supported by the National Natural Science Foundation of China
		(Grant No. 12031002 and No. 12671149). C. Zu was supported by the National Natural Science
		Foundation of China (Grant No. 12401151).

\subsection*{AI disclosure}

During the preparation of this manuscript, the authors
used ChatGPT (OpenAI) for English language polishing,
checking LaTeX syntax, and assisting with the
construction of the polynomial perturbations
\[
p(z)=v(z)+\frac{z^n}{\sqrt{w_n}}+\delta
\]
in Proposition~\ref{prop:perturbation}, including
the development of the associated uniform lower
bound estimates. The authors critically reviewed
and revised the AI-assisted mathematical arguments
and take full responsibility for the correctness
of all results and proofs presented in this paper.

	\medskip
	
	\noindent
	School of Mathematical Sciences, Dalian University of Technology,
	Dalian, Liaoning 116024, P. R. China
	
	\noindent
	Email address: \texttt{linzhaopeng2606@163.com}
	
	\medskip
	
	\noindent
	School of Mathematical Sciences, Dalian University of Technology,
	Dalian, Liaoning 116024, P. R. China
	
	\noindent
	Email address: \texttt{lyfdlut@dlut.edu.cn}
	
	\medskip
	
	\noindent
	School of Mathematical Sciences, Dalian University of Technology,
	Dalian, Liaoning 116024, P. R. China
	
	\noindent
	Email address: \texttt{zuchao@dlut.edu.cn}
	
\end{document}